\documentclass[final]{elsart}

\usepackage{amsmath}
\usepackage{amssymb,bm}
\usepackage{graphicx}
\usepackage{color}
\usepackage{mathrsfs} 
\makeatletter
\newif\if@restonecol
\makeatother

\usepackage{algorithm}
\usepackage{algpseudocode}
\DeclareMathAlphabet\mathpzc{OT1}{pzc}{m}{it}
\let\mathcal=\mathpzc
\def\E{{\mathbb E}}
\def\P{{\mathbb P}}

\def\B{\mathcal{B}}

\let\trueiiint=\iiint
\def\iiint{\mathop{\textstyle\trueiiint}\limits}
\def\intinfty{\int\limits_{\!\!-\infty\,\,}^{\,\,\infty\!\!}\kern-0.0em}
\def\iintinfty{\mathop{\int\!\!\int}\limits_{\!\!-\infty\,\,}^{\,\,\infty\!\!}\kern-0.0em}
\def\iiintinfty{\mathop{\int\!\!\int\!\!\int}\limits_{\!\!-\infty\,\,}^{\,\,\infty\!\!}\kern-0.0em}

\def\mse{{\mathrm{MSE}}}
\def\var{{\mathrm{VAR}}}

\let\<=\langle
\let\>=\rangle
\let\^=\hat
\def\~#1{{\mbox{\sf#1}}}
\let\'=\vec
\let\-=\mathbf

\def\circ{\ifmmode\mathchar"220E\else$\mathchar"220E$\fi}
\def\@#1{{\cal #1}}

\newcommand{\bx}{{\mathbf x}}

\numberwithin{equation}{section}

\journal{Elsevier}

\begin{document}
\centerline{}
\begin{frontmatter}



\title{A subset multicanonical Monte Carlo method for simulating rare failure events}


\author[authorlabel1]{Xinjuan Chen}
\ead{chenxinjuan@jmu.edu.cn}
\address[authorlabel1]{Department of Mathematics, College of Science, Jimei University,
Xiamen, Fujian,  361021, China.}

\author[authorlabel2]{Jinglai Li}
\ead{jinglaili@sjtu.edu.cn}
\address[authorlabel2]{Institute of Natural Sciences, Department of Mathematics, and MOE Key Laboratory of Scientific and Engineering Computing, Shanghai Jiao Tong University, Shanghai 200240, China. (Corresponding author)}


\medskip
\begin{center}
\end{center}

\begin{abstract}
Estimating failure probabilities of engineering systems is an important problem in many engineering fields.
In this work we consider such problems where the failure probability is extremely small (e.g $\leq10^{-10}$).  
In this case, standard Monte Carlo methods are not feasible due to the extraordinarily large number of samples required. 
To address these problems, we propose an algorithm that combines the main ideas of two very powerful failure probability estimation approaches: the subset simulation (SS) and the multicanonical Monte Carlo (MMC) methods. 
Unlike the standard MMC which samples in the entire domain of the input parameter in each iteration, 
the proposed subset MMC algorithm adaptively performs MMC simulations in a subset of the state space and thus improves the sampling efficiency.
With numerical examples we demonstrate that the proposed method is significantly more efficient than both of the SS and the MMC methods.
Moreover, the proposed algorithm can reconstruct the complete distribution function of the parameter of interest and thus can provide more information than just the failure probabilities of the systems. 
\end{abstract}


\end{frontmatter}

\section{Introduction}\label{s:intro}

Real-world engineering systems are unavoidably subject to various uncertainties
such as material properties, geometric parameters,
boundary conditions and applied loadings.
These uncertainties may cause undesired events, in particular, system failures or malfunctions, to occur.
Accurate evaluation of
failure probability of a given system is essential in many engineering fields such as risk
management~\cite{mcneil2015quantitative}, structural safety~\cite{robert1999structural},  {reliability-based} design and optimization~\cite{survey}, and thus is a central task of uncertainty quantification. 


Conventionally, the failure probability is often computed by constructing linear or quadratic expansions of the system model around
the so-called most probable point or $\beta$-point \cite{du2001most}, which is known as the first/second order reliability method~(FORM/SORM); see e.g., \cite{SchuellerPK04}
and the references therein.
It is well known that FORM/SORM may fail for systems with  {nonlinearity} or multiple failure modes.
The Monte Carlo (MC) simulation, which estimates the failure probability by repeatedly simulating the underlying system, provides
an accurate alternative to the FORM/SORM methods. The MC method does not make any reduction to the underlying system models, which means that it's applicable to any systems.
On the other hand, it is well known that the MC method suffers from slow convergence,
and thus can become prohibitively expensive
 when the system failures are rare (for example, around $10^{-10}$).
To this end, many advanced sampling schemes have been developed to reduce the estimation variance and improve the computational efficiency.
Among these schemes, the subset simulation (SS) method proposed by Au and Beck~\cite{AuB01,au2003subset}, is one of the most popular sampling strategies for estimating 
rare failure probabilities. 
Simply speaking, SS successively constructs a sequence of nested events with the very last one being the event of interest, 
and the probability of each event is estimated conditionally upon the previous one.
Other methods include,  just to name a few,  the cross entropy method~\cite{RubinsteinK_04,deBoerEtAl_AOR05,wang2015cross}, the population Monte Carlo~\cite{cappe2012population}.
Another attractive approach for estimating the failure probability is the  multicanonical Monte Carlo (MMC) method~\cite{berg1991multicanonical,berg1992multicanonical}, which was first developed to simulate rare events in physical systems. 
Later the method was used to estimate rare failure events in optical communication systems~\cite{holzloohner2003use,yevick2002multicanonical}. 
More recently, a surrogate accelerated  MMC method has been developed in \cite{wu2015surrogate} for uncertainty quantification applications.
The main idea of the MMC method is to partition the state space of the parameter of interest (which is usually a scalar and will be referred to 
 the performance parameter in what follows) into a set of small bins, and 
then iteratively construct a so-called flat-histogram distribution that can assign equal probabilities into each of the bins.  
Note that a major advantage of the MMC method is that it can reconstruct the entire distribution function of the parameter of interest,
and thus it can provide more information than just estimating the probability of a single event.

In this work, we propose a new algorithm that combines the key ideas of the SS and the MMC methods. 
Specifically, the new algorithm also constructs  a sequence of nested subdomains of the performance parameter, and then  performs the MMC scheme in each subdomain. The algorithm preserves some key properties of the standard MMC algorithm, while using
the subset idea to accelerate the computation. We thus refer to the proposed algorithm as the subset MMC (SMMC) method in the rest of the work.
Like the MMC method, the proposed SMMC algorithm can also compute the entire distribution function of the parameter of interest. 
Using several examples, we compare the performance of the proposed MMC algorithm with those of the SS and the MMC methods, and the numerical results show that SMMC method can significantly outperform both of the two original algorithms.

The rest of the work is organized as the following.
In Section~\ref{s:setup} we describe the mathematical formulation of the failure probability estimation problem. 
We then introduce the SS method in Section~\ref{s:ss} and the MMC method in Section~\ref{s:mmc} respectively. 
The proposed SMMC algorithm is presented in Section~\ref{s:smmc} and  three numerical examples are provided in Section~\ref{s:results}.  Some closing remarks will be given in Section~\ref{s:conclusion}.


\section{Failure probability estimation}\label{s:setup}

In this section, we shall describe the failure probability estimation problem in a general setting. Consider a probabilistic model where  $\-x$ is a $d$-dimensional random variable that represents the uncertainty in the model and the system failure is defined by a real-valued function  
\begin{equation}
y = f (\-x),
\label{perf}
\end{equation}
 which is known as the \emph{perform function}. 
For simplification, we shall assume that the state space of $\-x$ is $R^d$. 
The event of system failure is defined as that $y$ exceeds a certain threshold value $y^*$:   
\begin{equation}
F = \{\-x\in R^d\,|\, y=f (\-x) > y^*\},
\label{e:F}
\end{equation}
 and as a result the failure probability is
\begin{equation}
P_F = \P(F) = \int_{ \{\-x\in R^{d} | f(\-x)>y^*\}}  \pi(\-x) d\-x=\int_{\-x\in R^{d}} I_{F}(\-x) \pi(\-x) d\-x,\label{e:pf}
\end{equation}
where $I_{A}(\-x)$ is defined as an indicator function of set $A$:
\[
I_{A}(\-x) = \left\{ \begin{array}{ll}
         1 &\quad \mbox{if}\, \-x\in A,\\
         0 & \quad \mbox{if}\, \-x \notin A; \end{array} \right.
\]
and $\pi(\-x)$ is the probability density function (PDF) of $\-x$. In what follows we shall omit the integration domain when it is simply $R^{d}$. This is a general definition for failure probability, which is widely used in many disciplines involving with reliability analysis and risk management. Ideally, $P_F$ can be computed by using the standard MC estimation:
\begin{equation}
P_F\approx \frac{1}{N}\sum^N_{n=1}{I}_{\{f>y^*\}}(\-x_{n}),
\end{equation}
where samples $\-x_{1},...,\-x_{N}$ are drawn from the distribution with $\pi(\-x)$ as PDF.
However, as it has been discussed in Section~\ref{s:intro}, most engineering systems require high reliability, namely the failure probability $P_F \ll 1$. In this case, MC requires a large number of samples to produce a reliable estimate of $P_F$. 
On the other hand, in almost all practical cases, the performance function $f(\-x)$ does not admit analytical expression and has to be evaluated through expensive computer simulations, which makes the MC estimation of the failure probability prohibitive.Many advanced sampling schemes have been developed to compute the failure probability $P_F$, and we shall briefly introduce two popular choices of them: the SS and the MMC methods.

\section{The subset simulation method}\label{s:ss}
A brief introduction of the SS method, largely following \cite{AuB01}, will be provided in this section.
Note that we shall only outline the basic idea of the SS algorithm, and readers who are interested in the implementing details are referred to \cite{AuB01,au2007application,beck2015rare} and the references therein.

The idea of the SS method is to decompose the rare event $F$ into a sequence of  ``less-rare'' nested events,
\begin{equation*}
F = F_K \subset F_{K-1} \subset \cdots \subset F_1\subset F_0,
\end{equation*}
where $F_{k}$ is a more frequent event than $F_{k+1}$ for $k=1, \cdots, K-1$ and $F_0 = R^d$. Hence, the failure probability $P_F$ of the event $F$ can be computed by
\begin{eqnarray}
P_F=\mathbb{P}(F) = \mathbb{P}(F_K) &=& \mathbb{P}(F_1) \frac{\mathbb{P}(F_2)}{\mathbb{P}(F_1)} \frac{\mathbb{P}(F_3)}{\mathbb{P}(F_2)} \cdots \frac{\mathbb{P}(F_{K})}{\mathbb{P}(F_{K-1})}  \nonumber \\
&=& \mathbb{P}(F_1 |F_0) \mathbb{P}(F_2| F_1) \cdots \mathbb{P}(F_K | F_{K-1}),
\label{ss-condprod}
\end{eqnarray}
where $\mathbb{P}(F_k | F_{k-1})$ is the conditional probability of event $F_k$ given the occurrence of event $F_{k-1}$. Note that  $\mathbb{P}(F_1 |F_0)=\mathbb{P}(F_1)$. 

Before looking deeper into the algorithm, we will set up some new notations first.
Given an intermediate threshold value $y_{k}$,  we shall define $F_{k} = \{\-x\in R^d\,|\, f(\-x) > y_k\}$ as a corresponding intermediate event. In addition, we choose $y_0 = -\infty$ so that $F_0 = R^d$.  
The failure probability $P_F$ is now evaluated in a sequential manner.  In short words, starting from stage $k=0$,  
 the algorithm generates a  number of samples $\-x_{1}, \cdots, \-x_{N}$ from the distribution with PDF
\begin{equation}
\pi_k(\-x) =\pi \left( \-x | F_k\right) \propto {\pi(\-x) I_{F_k}(\-x)} ,
\label{ss-piF1}
\end{equation}
where it should be noted that $\pi_0(\-x) = \pi(\-x)$. 
It is worth noticing that drawing samples from $\pi_k(\cdot)$ is done with the Markov Chain Monte Carlo (MCMC) methods, which do not require the knowledge of the unavailable normalization constant in Eq.~\eqref{ss-piF1}. Afterward, one chooses an intermediate threshold value $y_{k+1}$ and  compute the conditional probability $\P(F_{k+1} | F_k)$ with standard MC, getting
\begin{equation}
\P(F_{k+1} | F_k) \approx \frac1N \sum_{n=1}^N I_{F_{k+1}}(\-x_{n}).
\end{equation}
The crucial point here is to choose the value of $y_{k+1}$ so that the resulting conditional probability $\P(F_{k+1} | F_k)$ is not too small. 
A commonly used approach is to let $y_{k+1}$ be the $(1-\gamma)$-th percentile of samples  $\{y_{1} = f(\-x_{1}), \cdots, y_{N} = f(\-x_{N})\}$ for some not too small positive number $\gamma$ (e.g., $=0.1$).
The algorithm proceeds until $y_{k+1}$ reaches $y^*$.
Therefore, one obtains the estimates of all the conditional probabilities $\P(F_{1}|F_0),\,\cdots,\,\P(F_K|F_{K-1})$ (assuming the algorithm reaches $y^*$ at the $(K-1)$-th iteration),  
and substituting the results into Eq.~\eqref{ss-condprod} yields an estimate of the desired failure probability $P_F$. 
We reinstate that the complete description of SS method is well documented in several works~\cite{AuB01,au2007application,Zuev2021,beck2015rare}.

\section{The multicanonical Monte Carlo method} \label{s:mmc}
We will now succinctly present the scheme of the MMC method, which is another effective algorithm used to estimate small failure probabilities.  Unlike the SS method, MMC solves the problem by constructing the distribution of the output parameter $y$. Namely, suppose that $\pi_y(\cdot)$ is the PDF of $y$, then the failure probability can be obtained by
\begin{equation}
P_F = \int_{y^*} ^{b} \pi_y(y) dy,
\end{equation}
where $b$ is in principle the maximum value of $y$. 
In practice, however, it is often not necessary to let $b$ be the maximum value of $y$, especially when $y$ is not bounded from above.   
It is easy to see that, for our purposes, it is sufficient to choose $b$ such that $\P(y>b)\ll\P(y>y^*)$. 
Hence, in order to find the failure probability of the system, one only needs the PDF of $y$. To be more precise, we only need the PDF of $y$ in the interval $[y^*,b]$.
This is not a simple task, however, because the failure region is typically located in the tail of $y$. 
 
A popular strategy applied to estimate the PDF of a continuous random variable $y$ with simulation  
is to approximate the PDF with histograms. 
Suppose we are interested in the PDF of $y$ in the interval $\B=[a,\,b]$, and we first equally decompose $\B$ into $m$ bins of width $\Delta$, whose centers are the discrete values $\{b_1, ... , b_m\}$. We define the $i$-th bin as the interval $B_i = [b_i  - \Delta/2, b_i + \Delta//2]$ and the probability for $y$ to be in $B_i$ is $P_i =\P\{y \in B_i\}$. 
Note that the width of each bin needs not to be identical in principle, and here we use identical bin width just for the simplicity of notations. 
The PDF of $y$ at point $b_i$ then can be approximated by \[p(b_i) \approx P_i / \Delta,\] if $\Delta$ is sufficiently small. This binning implicitly defines a partition of the input space $X$ into $m$ domains $\{ D_i \}^{m}_{i = 1}$, where
$$D_i = \{ \-x \in  R^d: f(\-x) \in B_i \}$$
is the domain in $X$ that is mapped into the $i$-th bin $B_i$ by $f(\-x)$. Note that, while $B_i$ are simple intervals, the domains $D_i$ are multidimensional regions with possibly tortuous topologies. 
As a result, the probability $P_i$ can be re-written as an integral in the input space:
\begin{equation}
P_i = \int_{D_i} \pi (\-x) dx = \int I_{D_i}(\-x) \pi (\-x) dx = \E[I_{D_i}(\-x)].\label{eq:w:1}
\end{equation}
Now suppose that $N$ samples $\{\-x_1,\ldots,\-x_N\}$ are drawn from the distribution $\pi(\-x)$, possibly with MCMC, 
$P_i$ can be estimated with the MC estimator:
\begin{equation}
\hat{P}^{MC}_i = \frac{1}{N} \sum^{N}_{j = 1} I_{D_i}(\-x_j) = \frac{N_i}{N}, \label{eq:w:2}
\end{equation}
where $N_i$ is the number of samples that fall in the domain $D_i$. 

As it's well known, standard MC simulations have difficulty in reliably estimating the probabilities in the tail bins. The technique of importance sampling~(IS) can be effectively used to address the issue. The principle idea of IS is to choose a biasing distribution $q(\-x)$ and rewrite Eq. \eqref{eq:w:1} as 
\begin{equation}
P_i = \int I_{D_i}(\-x) [\frac{\pi(\-x)}{q(\-x)}] q (\-x) d\-x = \E^{*} [I_{D_i}(X) w(X)]
 \label{eq:w:3}
\end{equation}
where  $w(\-x)  = \pi(\-x)/q(\-x)$ is called the IS weight, and $\E^*$ indicates expectation with respect to the biasing distribution $q(\-x)$. It follows that the IS estimator of $P_i$ becomes
\begin{equation}
  \hat{P}^{IS}_i  =  {\left( \frac{N^*_i}{N} \right)}  {\left[ \frac{1}{N^{*}_i} \sum^{N}_{j = 1} I_{D_i}(\-x_j)w(\-x_{j})\right]},
   \label{eq:w:5}
\end{equation}
where the samples $\{\-x_1,\ldots,\-x_N\}$ are now generated from the biasing distribution $q(\-x)$,
and $N^*_i$ is the number of samples falling in the region $D_i$. 

One can easily see that the key of IS is to choose an appropriate biasing distribution $q(\-x)$ that can help to achieve the objective of the simulation. While regular IS usually aims to estimate the probability in a given region,
the goal of our simulation is to have a good estimate of $P_i$ for all $i=1\ldots m$, and in this respect, it is reasonable to seek a biasing distribution that assigns equal probability to each bin and zero probability for any region outside $\@D = \cup_{i=1}^m D_i$, which implies that
\begin{subequations}
\label{e:pis}
\begin{equation}
P^*_1=P_2^*=...P^*_m=1/m,
\end{equation}
where
\begin{equation}
P^{*}_i = \int_{R^d} I_{D_i} (\-x) q(\-x) d\-x = E^{*} [I_{D_i} (X)], \mbox{for}\  i=1, \cdots, m.
\end{equation}
\end{subequations}

We refer to the biasing distribution which satisfies Eqs.~\eqref{e:pis} to be flat-histogram (FH). One should be noted that the FH distributions  are not unique, and among them there is one which assigns a constant weight to all $\-x \in D_i$, i.e. $w(\-x) = w_i$ for $\-x\in D_i$ where $w_i = P_i / P^*_i$. In this case, the biasing distribution $q(\-x)$ is called to be uniform-weight (UW).

In particular, we assume that the biasing distribution $q(\-x)$ is given in the form of

\begin{align}
 \label{eq:w:6}
q(\-x) = \left\{
\begin{aligned}
 & \frac{\pi(\-x)}{c_\Theta \Theta(\-x)}   & \-x \in \@D;  \\
&0 & \-x\notin \@D,
\end{aligned}
\right.
\end{align}
where $\Theta(\-x) = \Theta_i>0$ for all $x \in D_i, i = 1,...,m$, satisfying 
\begin{equation}
\sum_{i=1}^n \Theta_i = 1, \label{e:sumTheta}
\end{equation}
 and $c_\Theta$ being a normalized constant.
 
 It is easy to show that $q(\-x)$ given in Eq.~\eqref{eq:w:6} is UW with  $w_i = c_\Theta \Theta_i$ for $i=1,\ldots, m$. Next we shall impose the constraint so that $q(\-x)$ given in Eq.~\eqref{eq:w:6} is FH. Since 
\begin{equation}
P^*_i = \int_{D_i} q(\-x)d\-x = \frac{\int_{D_i} \pi(\-x) d\-x}{c_\Theta \Theta_i} = \frac{P_i}{c_\Theta \Theta_i} ,  
 \label{eq:w:7}
\end{equation}
and by setting the left hand side of Eq.~\eqref{eq:w:7} to be qual to $1/m$, we obtain 
\begin{equation}
\Theta_i = \frac{m}{c_\Theta} P_i. \label{e:theta_i}
\end{equation}
Substituting Eq.~\eqref{e:theta_i} into Eq.~\eqref{e:sumTheta} results in $c_\Theta = m \rho$, where $\rho = \sum_{i=1}^m P_i$, and it follows immediately that $\Theta_i = P_i/\rho$. Note that in general the probability $\rho = \P[y\in \B]\leq1$ and is unknown in advance. 
A conventional solution is to take a sufficiently large interval $\B$ so that  $\rho\approx1$, and we adopt this choice in this work.

\begin{figure}
\centerline{\includegraphics[width=.9\textwidth]{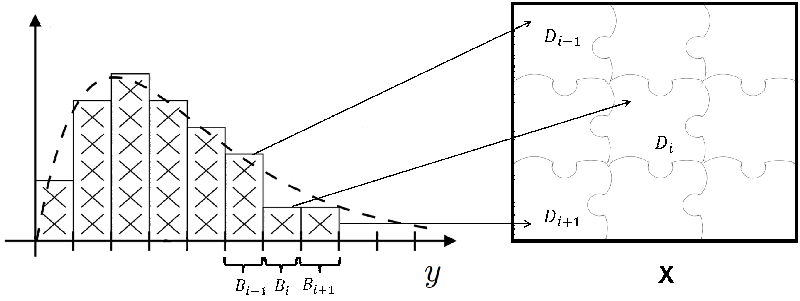}} \label{fig:bins}
\caption{Schematic illustration of the connection between $B_i$ and $D_i$.}
\end{figure}

However, for the reason that $\Theta_i$, $i = 1, \cdots, m,$ depend on the sought after unknown $P_i$,   the actual UW-FH distribution just derived above can not be utilized directly to achieve the goal of sampling equally in each bin. 

The MMC method uses an adaptive scheme to address this issue. Simply speaking, MMC adaptively constructs a sequence of distributions 
\begin{align}
 \label{e:qn}
q_k(\-x) = \left\{
\begin{aligned}
 &\frac{\pi(\-x)}{c_k \Theta_k(\-x)},    & \-x \in \@D;  \\
&0 & \-x\notin \@D, 
\end{aligned}
\right.
\end{align}
where $\Theta_k(\-x) = \Theta_{k,i}$ for $\-x \in D_i$, converging to the actual UW-FH distribution.
Before proceeding to the MMC algorithm, we derive an alternative representation of $\Theta_i$ from Eq. \eqref{e:theta_i}: 
\begin{equation}
\Theta_i =  P^*_i w_i/\rho, \quad \mathrm{for} \quad i=1,\cdots, m. \label{e:theta_i2}
\end{equation}
Typically, the MMC method starts from the original PDF $q_0(\-x) = \pi(\-x)$, where the associated parameter values are $c_0=1$ and $\Theta_{0,i} = 1$ for all $i=1,\ldots,m$. 
In the $k$-th iteration, one first draws $N$ samples $\{ \-x_j \} ^N_{j = 1}$ from the current distribution $q_k(\-x)$, 
and then updates $\{\Theta_{k+1,i}\}_{i=1}^m$ using the following formulas derived from Eq.~\eqref{e:theta_i2}:
\begin{subequations}
\label{e:params}
\begin{gather}
\hat{H}_{k,i} = \frac{N^*_{k,i}}{N},\\
w_{k,i} = c_k\Theta_{k,i},\\
\Theta_{k+1,i} 
= \hat{H}_{k,i}  w_{k,i}/\rho,  
\label{eq:w:8}
\end{gather}
\end{subequations}
where $N^*_{k,i}$ is the number of samples falling into region $D_i$ in the $k$-th iteration. 
It should be noted that, MMC usually employs MCMC to draw samples from $q_k(\-x)$, thanks to which we do not need to estimate $c_k$ during the iterations (i.e., just to take $c_k=1$ in each iteration). 
However, the constant is needed when one wants to compute $P_i$ for $i=1, \cdots, m$, using the IS estimator~\eqref{eq:w:5}, in the final stage.  To circumvent the obstacle, we estimate $P_i$ by
 \[P_i \approx \frac{\Theta_{K,i}}{\sum_{i=1}^m{\Theta_{K,i}}} \rho ,\quad \mathrm{for}\quad i=1, \cdots, m,\]  
where $K$ is the index of the final iteration.
Formal convergence analysis, as well as possible improvements of the MMC method are not discussed in this work, and readers who are interested may consult, e.g.\cite{berg2000introduction,berg2004markov,iba2014multicanonical,landau2014guide}, and the references therein.

\section{The subset MMC method} \label{s:smmc}
As it has been described in the previous section, the conventional MMC method uses a sufficiently large interval $\B$ such that $\rho=1$, which, unfortunately, is not an efficient approach for our purposes because the failure region that we are interested in, $[y^*,\,b]$, is typically a small subinterval of $\B$.  As a result, only a very small portion of the samples will be used to estimate the density in the region of interest. To address the issue, we propose a subset MMC algorithm,  which combines the key ideas of the SS and the MMC methods.

Consider the case where $\rho$ is unknown.  Similar to the SS method, we now construct a sequence of nested intervals $\B^0 \supset...\supset \B^J$, where $\B^0=\B$ and $\B^J =[y^*,\,b]$ is the interval of interest.  It should be obvious to see that the corresponding domains in the input space are also nested. Let $\rho_j = \P(y\in \B^j)$ for $j = 1, \cdots, J$, and as  it's  explained before, the failure probability $P_F\approx\rho_J$.

Let the bins $B_1,...,B_m$ be predetermined as the previous section and will not be changed as the algorithm proceeds. Moreover, for the sake of simplicity, we assume that the threshold value $y^*$ coincides with the left boundary of one of the bins, i.e. $b_{m^*} - \Delta/2 =y^*$ for some integer $1\leq m^*\leq m$.  In this case, it will be natural to construct each interval $\B^j$ as a union of bins: $\B^j = \cup_{i={m_j}}^m B_i$ for some integer $1\leq m_j\leq m^*$. It can be easily seen that $m_1\leq m_2\leq \cdots $. Starting from $\B^0$ (with $\rho_0=1$ and $m_0=1$), we now perform a standard MMC within the interval $\B^j$ and compute the probabilities of bins from $m_j$ to $m$: $P_{m_j},...P_{m}$.   An $m_{j+1}$ is chosen such that $m_j\leq m_{j+1}\leq m$, which indicates that the choice of $m_{j+1}$ determines the next interval $\B^{j+1}$ (the criterion that we use to determine $m_j$ will be provided later).  The basic thought is that we can gradually concentrate the samples toward the region of interest. It follows immediately that the associated probability $\rho_{j+1}$ can be estimated by
\[ \rho_{j+1} = \sum_{i={m_{j+1}}}^m P_i.\]
The algorithm proceeds until $m_{j+1} =m^*$. The complete scheme is described in Algorithm~\ref{smmc}. Several remarks regarding the implementation of the algorithm of SMMC are given in the following. 
\begin{itemize}
\item In each iteration, $m_{j+1}$ is determined according to the following. Firstly,  a not-too-small positive number $\alpha<1$ (e.g. $=0.2$) is chosen. 
Then we select an $m_{j+1}$ such that (approximately) $100\alpha\%$ of the samples fall in the interval $\@B^{j} = \cup_{i={m_{j+1}}}^m B_i$.
\item In line~\ref{ln:mcmc}, the samples are drawn from $q_k(\cdot)$ using the MCMC methods. In particular, we implement a multiple chain MCMC algorithm specifically tailored for this problem.
The details of the algorithm is given in Appendix~A.  
\item The terminating condition used here is $m_j=m^*$, i.e., when the interval in which we perform MMC reaches the area of interest. 
\end{itemize}

\begin{algorithm}
\caption{The subset MMC algorithm}\label{smmc} 
\begin{algorithmic}[1]
\Require{$\pi({\bf x})$, $\{ B_i\}_{i=1}^m$, $m^*$, $n$, $K$, $\alpha$.}
\Ensure{$P_F$.}
\Procedure{$P_F=$ SMMC}{$\pi(\bx)$, $\{ B_i\}_{i=1}^m$, $m^*$, $n$, $\alpha$}
\State Initialization: $j=0$, $\bm\Theta^0(\bx) = \left(1, \cdots, 1 \right)$, $m_0=1$, $\rho_0 = 1$; 
 \While{$m_j< m^*$}
\State $[\bm\Theta^{j+1},m_{j+1}] = \mathrm{MMC}(\bm\Theta^j, m_j, \pi(\bx), \{ B_i\}_{i={m_j}}^m, m^*, n, K, \alpha)$;
\For{$i=m_{j}...m$}
\State $P_i = \frac{\Theta^{j+1}_{i}}{\sum_{i={m_j}}^m{\Theta^{j+1}_{i}}} \rho_{j+1}$;
\EndFor
\State  $\rho_{j+1} = \sum_{i={m_{j+1}}}^m P_i$; 
\State $j=j+1$;

\EndWhile
\State $P_F = \rho_j$;  
\EndProcedure
\Procedure{$[\bm\Theta^+, m^+]=$MMC} {$\bm\Theta^-,m^-, \pi(\bx),\{ B_i\}_{i=1}^m, m^*, n, K, \alpha$}
\State $\bm\Theta_0 = \bm\Theta^-$;
\State $D = \{\bx\in R^d\,|\, g(\bx) \in \cup_{i={m^-}}^m B_i\}$;
\For{$k=0...K$} 
\State Let $q_k$ be given by Eq.~\eqref{e:qn} with $\bm\Theta_k$;
\State Draw $n$ samples $\{\bx_1,...,\bx_n\}$ from $q_k$; \label{ln:mcmc}
\State Evaluate $S_k = \{f(\bx_1),...,f(\bx_n)\}$;
\State Compute $\bm\Theta_{k+1}$ using Eqs.~\eqref{e:params};
\EndFor
\State Let $y_\alpha$ be the $(1-\alpha)$-th quantile of set $S_K$;
\State Let $m^+$ be the index of the bin such that $y_\alpha\in B_{m^+}$;
\State $m^+ = \min\{m^+,\,m^*\}$;
\For{$i=m^+...m$}
\State $\Theta^+_{1+i-m^+} = \Theta^{K+1}_{i}$;
\EndFor
\EndProcedure
\end{algorithmic}

\end{algorithm}
So far, although we have presented the SMMC algorithm as a variant of the standard MMC method, 
the key idea is also inspired by the SS method, given that, in the SMMC algorithm, a sequence of subsets are first constructed and in each subset a MMC iteration rather than plain MC is performed in order to drive samples towards the failure region. It is also interesting to notice that, just like the  standard MMC method, one of the significant advantages of the SMMC method is that, if desired, it can construct the entire PDF of $y$ without any additional cost, as we obtain the estimates of $P_i$ for $i=1, \cdots, m$ during the iteration.
We shall illustrate this advantage with numerical examples in Section~\ref{s:results}.

Finally we shall provide a simple analysis of the estimator error of the SMMC algorithm. 
Two simplifications are made for convenience. One is that the samples drawn are independent while noting that the samples are certainly not independent when they are drawn with MCMC methods.
The other assumption is that  the biasing distribution is ``perfectly flat'' in the last iteration, namely,  
the biasing distribution is  given by Eq.~\eqref{eq:w:6}, where $\{\Theta_i\}_{i={m^*}}^m$ are given by Eq.~\eqref{e:theta_i}
and 
\[c_\Theta=(m-m^*)\rho,\]
where $(m-m^*)$ is the number of bins in the last iteration.
  Now for any $m^*\leq i\leq m$, the estimator of $P_i$ is 
\begin{equation}
\hat{P}_i = \sum_{j=1}^N I_{D_i}(\-x_j) w(\-x_j), \label{e:Phat_i}
\end{equation}
where the samples are drawn from distribution Eq.~\eqref{eq:w:6} with
\begin{equation}
w (\-x)= (m-m^*) \rho' \Theta(\-x).\label{e:weight}
\end{equation}
Note that $\rho'$ in Eq.~\eqref{e:weight} is an estimate of $\rho$ as 
the actual value of $\rho$ is unknown in our problem.

The mean square error (MSE) of Eq.~\eqref{e:Phat_i} is computed as 
\begin{eqnarray*}
\mse[\hat{P}_j] &= &\var[\hat{P}_j]+(\E[\hat{P_j}]-P_j)^2,\\
&=& \frac1n \left((\frac{\rho'}\rho)^2 (m-m^*) P^2_j-(\frac{\rho'}\rho)^2 P_j^2\right)+ \left((\frac{\rho'}\rho) -1\right)^2P_j^2,\\
&=&  \frac{(m-m^*-1)}N \phi^2P_j^2 +(\phi-1)^2P_j^2 ,
\end{eqnarray*}
where $\phi = {\rho'}/\rho$.
It is not difficult to see that the optimal value of $\phi$ that minimizes the MSE is $\phi = N/(N+m-m^*-1)$, and the resulting minimal MSE is 
\[\mse_{\min} = \frac{m-m^*-1}{m-m^*-1+N} P^2_j.\]
It is interesting to see from the results that for the MSE to be minimal, one should choose 
\[{\rho'} = \frac{N}{(N+m-m^*-1)}\rho,\]\
rather than ${\rho'}=\rho$. 
However, when $m-m^*\ll N$ which is the usual case, 
$\frac{N}{(N+m-m^*-1)} \approx 1$, and thus we  choose not to include the factor $\frac{N}{(N+m-m^*-1)}$ in the estimate of 
$\rho$ in Algorithm~1. 

\section{Numerical examples} \label{s:results}
 
\subsection{A two-dimensional mathematical example}
The first example is a two-dimensional mathematical problem. 
Suppose that $\-x= (x_1, x_2)$ is a two-dimensional random variable where $x_1$ and $x_2$ both follow standard normal distribution and are independent to each other.  The event of failure is defined as 
\[ \min\{\|\-x-\-x_r\|, \|\-x-\-x_l\|\}<1,
\]
where $\-x_r=(8,2)$ and $\-x_l=(-8,2)$.
On the one hand, it is clear that the problem has two disjoint failure domains: $\{x\in R^2 | \|\-x-\-x_r\|<1\}$ and $\{\-x\in R^2|\|\-x-\-x_l\|<1\}$,  
which poses challenge for many standard IS methods. 
On the other hand, for this two-dimensional example, the failure probability can be accurately estimated by performing a numerical integration, yielding $P_F = 1.41\times10^{-13}$. 

We estimate the probability with three methods: SS, standard MMC and the proposed SMMC algorithms. As for the SS  method, we largely follow the implementation described in \cite{Zuev2021} and set $\gamma=0.1$ as it's suggested in \cite{Zuev2021}. For both of the standard MMC and proposed SMMC methods,  the entire region of interest of the output is taken to be $[0,100]$ and is equally divided into 100 bins. For the SMMC method, we take $\alpha=0.2$ to generate the nested intervals. 

\begin{figure}
\centerline{\includegraphics[width=1\textwidth]{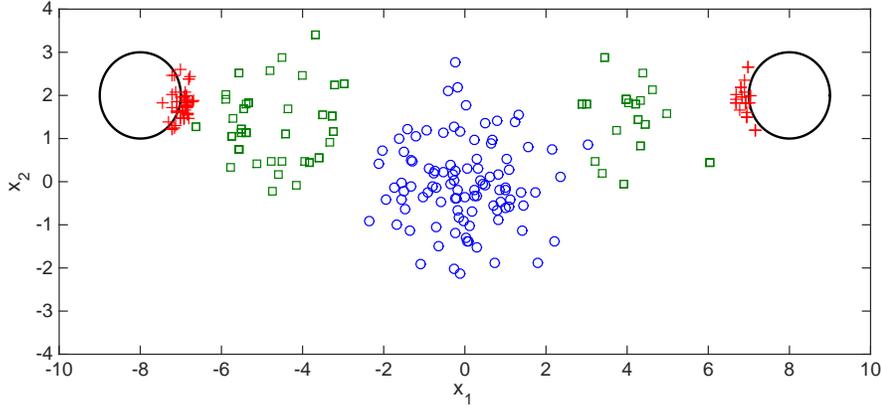}} \label{fig:2D}
\caption{Schematic illustration of the connection between $B_i$ and $D_i$.}
\end{figure}
We compare the three methods using three different numbers of points: $2\times 10^4$, $1\times10^5$ and $3\times10^5$. 
It should be noted that, while for the MMC algorithm we can choose the exact sample size by fixing 
the total number of iterations and the number of samples used in each iteration,  we can not exactly control  it in the SS and the SMMC methods, 
and so we can only adjust the algorithms so that the total amounts of samples are close to the aforementioned numbers.  Also, during all the computations, we manage to ensure that the sample size of the SMMC is smaller than those in the other two algorithms. For each of the three methods, we repeatedly perform the simulations for 100 times and computed the average number of samples 
as well as the  relative mean square errors (RMSE):
\[ \mathrm{RMSE} =\frac{ \frac1L\sum_{l=1}^L|\hat{P}_l-P_F|^2}{P_F^2},
\]
where $L=100$ is the total number of computations and $\hat{P}_l$ is the estimated probability at the $l$-th test. 
The test results are presented in Table \ref{tb:P1results}. We have found that, in the SMMC method, most tests terminate within three iterations. We also show the sample distributions in one test trial in Fig~\ref{fig:2D}. One can see from the figure that, the SMMC method is capable of directing samples toward the failure region in a rather efficient manner.  More information can be learned from the results in Table \ref{tb:P1results}.  In particular, we can see that in all the cases, the SMMC performs substantially better than the other two methods, even with less samples.

\begin{table}[!htb]
\center
\ {
\begin{tabular}{lccc}
\hline %
&SS& MMC&   SMMC\\
\hline %
$n$ & $2.54\times10^4$ & $2.0\times10^4$ & $1.70\times10^4$ \\ \cline{2-4}  
RMSE & 49  & 25.9 & 0.49 \\ \cline {1 -4} 
 $n$ & $1.14\times10^5$ & $1.0\times10^5$ & $0.98\times10^5$\\ \cline{2-4}  
RMSE & 2.69  &  2.56 & 0.078 \\ \cline {1 -4} 
$n$ & $3.34\times10^5$ & $3.0\times10^5$ & $2.99\times10^5$\\ \cline{2-4}  
RMSE & 0.16  &  2.13 & 0.017 \\ \cline {1 -4} 
  \hline
 \end{tabular}
}
\caption{Example 1: performance comparison of the three methods with different sample sizes.}
\label{tb:P1results}
\end{table}

\subsection{A high dimensional mathematical example}
This one is also a mathematical problem, but of a higher dimensionality than the previous one.  Specifically, we let $\-x$ be a $d$-dimensional random variable following 
standard Gaussian distribution: $\-x \sim N(0,I)$ where $I$ is the $d\times d$ identity matrix. The failure event is defined as
$f(\-x)>y^*$ with 
\begin{equation}
f(\-x) =  \|\-x\|^2_2.
\end{equation}
In our numerical tests, we choose $d=10$ and $y^*=75$. In this setting, the failure probability  can be computed analytically as $P_F = 4.76\times10^{-12}$. The challenge of this example is that, it is rather difficult to construct an effective parameterized form for the biasing distribution, which is a critical issue for most of IS methods.  A biasing distribution in a Gaussian form may not perform well for this problem, just to name a few.

\begin{figure}
\centerline{\includegraphics[width=0.9\textwidth]{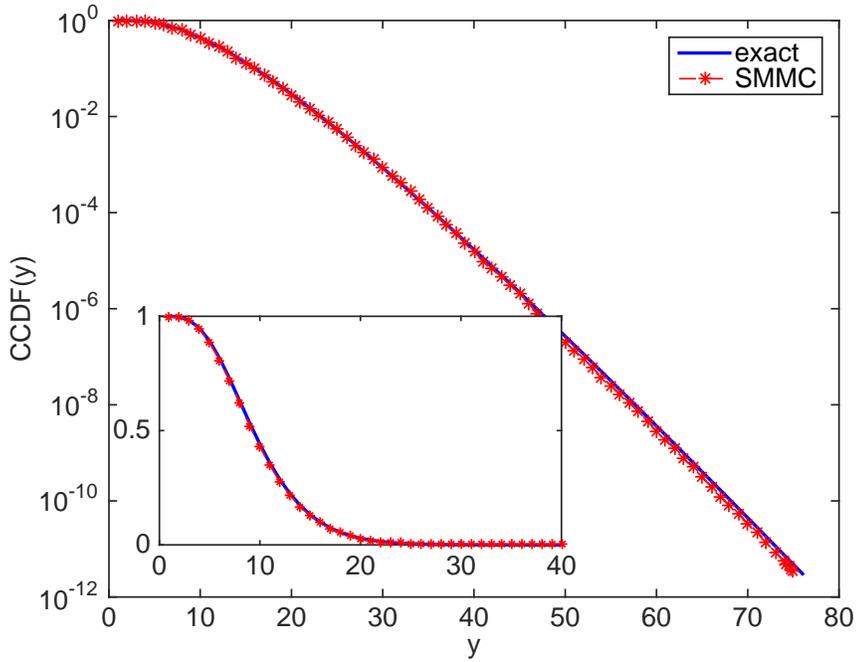}} \label{fig:ccdf}
\caption{The CCDF computed by the SMMC method compared to the exact results, both are plotted on a logarithmic scale.
Inset: the same plots but on a linear scale.}
\end{figure}

We tested all the three methods on this example as well. The specifications of the implementations of all the three methods are kept the same as those in the first example. Like the first example,  we test each method with three different sample size: $2\times10^4$, $1\times10^5$ and $3\times10^5$, 
and  repeatedly perform the simulations 100 times for each sample size. 
The RMSE results are available in Table~\ref{tb:P2results}. 
The results indicate that,  in this example, the SMMC method also substantially outperformed the other two methods. Besides, as it's mentioned earlier, another improvement of the SMMC method over SS method is that it can also be used to construct the complete distribution of the output $y$. To show this, we plot in Fig.~\ref{fig:ccdf} the complement cumulative distribution function (CCDF) of $y$ obtained by the SMMC method,  which is defined as 
\[\mathrm{CCDF}(y) = 1 - \mathrm{CDF}(y),\]
where $\mathrm{CDF}(y)$ is the cumulative distribution function (CDF) of $y$.
As a comparison, we also show the exact CCDF function of $y$, and one can see that the result of SMMC agrees very well with the exact one.

\begin{table}[!htb]
\center
\ {
\begin{tabular}{lccc}
\hline %
&SS& MMC&   SMMC\\
\hline %
$n$ & $2.5\times10^4$ & $2.0\times10^4$ & $1.7\times10^4$ \\ \cline{2-4}  
RMSE & 6.3  & 7.1 & 0.5 \\ \cline {1 -4} 
 $n$& $1.09\times10^5$ & $1.0\times10^5$ & $1.07\times10^5$\\ \cline{2-4}  
RMSE& 2.6  &  2.7 & 0.2 \\ \cline {1 -4} 
 $n$& $3.27\times10^5$ & $3.0\times10^5$ & $2.90\times10^5$\\ \cline{2-4}  
RMSE & 0.15  &  2.4 & 0.02 \\ \cline {1 -4} 
  \hline
 \end{tabular}
}
\caption{Example 2: performance comparison of the three methods with different sample sizes.}
\label{tb:P2results}
\end{table}

\subsection{Quarter car model}
\begin{figure}
\centerline{\includegraphics[width=0.5\textwidth]{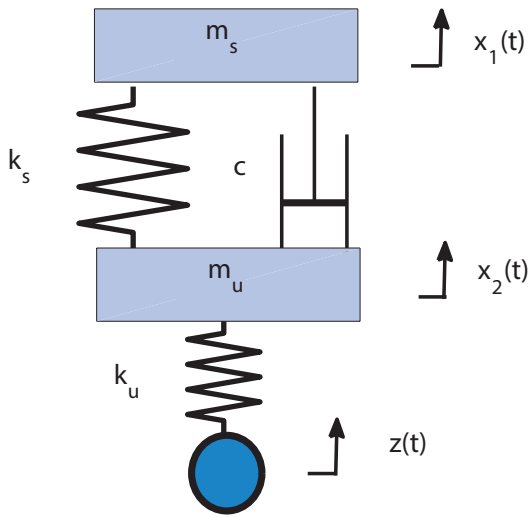}} \label{fig:qcar}
\caption{The schematic illustration of the  quarter car model. }
\end{figure}

The last example is the quarter car model for vehicle suspension systems~\cite{wong2001theory}.
The schematic illustration of the model is shown in Fig~\ref{fig:qcar}, where the sprung mass $m_s$ and the unsprung mass $m_u$ are connected by a nonlinear spring and a linear damper.
The stiffness of the nonlinear spring is $k_s$ and the damping coefficient of the linear damper is $c$. 
The displacement of the wheel $z(t)$ represents the interaction of 
the  quarter  car  system  with  the  terrain.
Mathematically, the model is described by a  two-degree-of-freedom ODE system~\cite{wong2001theory}:
\begin{subequations}
\label{e:qcar}
\begin{align}
 m_s \frac{d^2x_1}{dt^2} = & - k_s (x_1-x_2)^3 - c(\frac{dx_1}{dt}-\frac{dx_2}{dt}),\\
m_u\frac{d^2x_2}{dt^2} = & ~k_s (x_1-x_2)^3  
+c(\frac{dx_1}{dt}-\frac{dx_2}{dt})+k_u(z(t)-x_2).
\end{align}
\end{subequations}
where $x_1$ and $x_2$ are the displacements of the sprung and the unsprung masses respectively. 
In our example, we assume that the uncertainty in the system arises from the random road profile, and as a result the wheel displacement $z(t)$ is modeled 
as a zero-mean white Gaussian random force with standard deviation $\sigma = 0.05$. 
The other model parameters are all taken to be fixed and the values of them are shown in Table~\ref{tb:P3params}. 
The quantity of interest is the maximum difference between displacements of the sprung and the unsprung springs in a given interval $[0,\,T]$,  
\[y = \max_{0\leq t\leq T}\{|x_1(t)-x_2(t)\}|,\]
and we want to reconstruct the CCDF of $y$. 
With the CCDF, we can estimate directly the probability $\P(y>y^*)$ for any $y^*$ in the range of interest. 

In the numerical simulations, we take $T=1$, and the initial conditions of Eqs.~\eqref{e:qcar} to be
\[
x_1(0) = \frac{dx_1}{dt}(0)=0,\quad
x_2(0) = \frac{dx_2}{dt}(0)=0.
\]
The Eqs.~\eqref{e:qcar} is numerically solved with the classical Runge-Kutta method where the step size is taken to be $\Delta t = T/100$, which means that the random variable in this problem is effectively of 100 dimensions. 

\begin{figure}
\centerline{\includegraphics[width=0.9\textwidth]{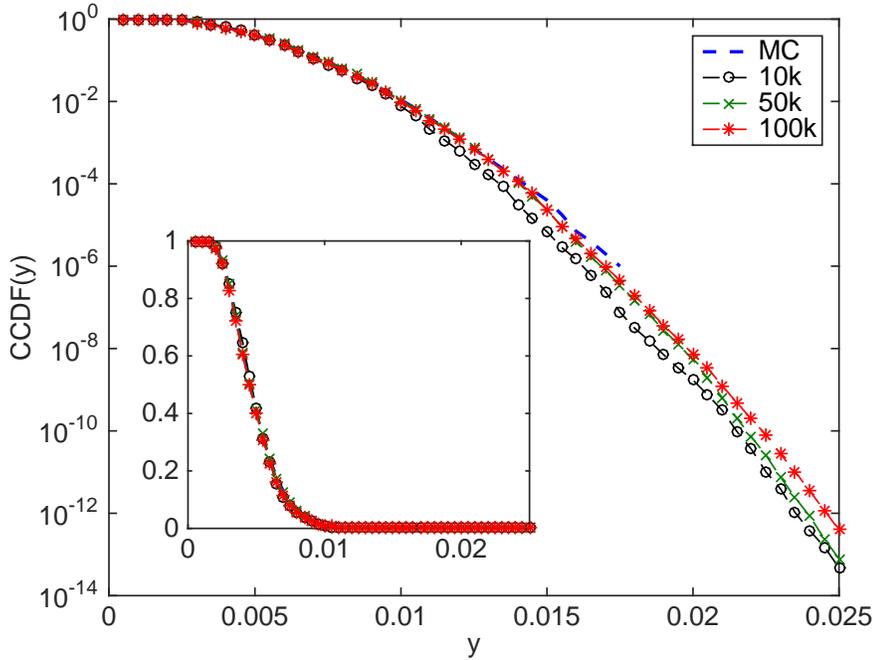}} \label{fig:ccdf_qcar}
\caption{The CCDF computed by the SMMC method with three different sample sizes: $10^4$ (circles), $5\times10^4$ (crosses) and $10^5$ (asterisks).
As a comparison, we also plot in the figure the result of standard MC with $10^6$ samples (dashed line). All are plotted on a logarithmic scale.
Inset: the same plots but on a linear scale. }
\end{figure}

Also, the CCDF of $y$ using a standard MC method with $10^6$ samples is constructed. 
We perform the SMMC method with three sample sizes $10^4$, $5\times 10^4$ and $10^5$ respectively, and present all the results in Fig.~\ref{fig:ccdf_qcar}.
One can see from the figure that, the results of the SMMC agree largely with those of the standard MC. Without surprise, the MC method can only obtain the CCDF at the order of $10^{-6}$, while the SMMC method can compute the CCDF  down to $10^{-12}$ and smaller with much less samples than the MC method.  One can also see that the result of the SMMC of $10^4$ samples departs evidently from those of $5\times10^4$ and $10^5$, indicating that  the sample size of $10^4$ may not be sufficient for this problem. With around $10^5$ samples, we can compute the probability as small as $10^{-12}$ using the SMMC method. Note that the CCDF computed with the SMMC method can also provides us with other important information such as the extreme quantiles.
For instance, we can see directly from the CCDF that the $(1-10^{-8})$-th quantile is 0.0198 and 
the $(1-10^{-10})$-th is 0.0224.  Such information can not be easily obtained with the SS method. 

\begin{table}[!htb]
\center
\ {
\begin{tabular}{lcccccc}
\hline %
&$m_s$& $m_u$&  $k_s$ & $k_u$ & $c$ \\
\hline %
 & $20$ & $40$ &400 & 2000 & 600\\
  \hline
 \end{tabular}
}
\caption{The parameter values of the quarter car model.}
\label{tb:P3params}
\end{table}


\section{Conclusions}\label{s:conclusion}
In summary, we propose an efficient algorithm for estimating failure probabilities of complex engineering systems, which combines the central ideas of the SS and the MMC methods.
The new algorithm constructs a sequence of subdomains of the performance parameter $y$ and  performs regular MMC iterations within each subdomain only. We demonstrate that the proposed SMMC method can significantly outperform the two original methods, and moreover, like the MMC method, it can be used to reconstruct the entire distribution function of the performance parameter. 
We believe that the SMMC method can be a useful tool for many practical engineering problems that involve failure probability estimations. 

Several improvements and extensions of the proposed algorithm are possible. 
First, for systems with highly intensive computer models, even with the SMMC method, the total computational cost is still unaffordable. In such problems, a possible solution is to construct computationally inexpensive surrogate models and use them in  the simulations (see, e.g. \cite{li2011efficient,LiXiu_JCP10,dubourg2013metamodel}). To this end, surrogates have been used to accelerate the simulations in both the SS \cite{papadopoulos2012accelerated} and the MMC~\cite{wu2015surrogate} methods. Thus we hope to develop surrogate based methods to reduce the computational cost of the SMMC algorithm. 
Secondly, in many practical problems, we often have  computer models with different fidelities for the system. In this case, a very interesting question will be how to incorporate the multi-fidelity models  with the SMMC algorithm and further improve the computational efficiency. 
Finally we think the proposed method can also be applied to problems beyond failure probability estimations. In particular, we hope to apply the  SMMC algorithm with necessary modifications to evaluate the evidence (normalization constant of the posterior distribution) in Bayesian inference problems. We plan to study these problems in future works. 

\section*{Acknowledgment}
The work was partially supported by the National Natural Science Foundation of China under grant number 11301337.

\appendix

\section{A specialized MCMC algorithm for the SMMC simulations}

Here we present a specialized  MCMC algorithm, largely following the modified Metropolis algorithm used in the SS method~\cite{AuB01}. 
First, unlike the standard MMC algorithm which employs only one MCMC chain  at each cycle, we uses a multi-chain MCMC algorithm. 
In particular, in each cycle we randomly select a sample from each $D_{1}, \,D_{2},\, \cdots,\, D_m$ if there are any,
and then we use the obtained $m' \leq m$ samples as the seeds to perform $m'$ chains parallely. 
Note that here $m'$ is automatically determined by the algorithm and for this reason, we can not strictly specify the number of samples drawn
in each MMC iteration. 

Next we adopt the dimension by dimension proposal used in \cite{AuB01}. 
To do so, we need to assume that in the original distribution $\pi(\-x)$ all the component of $\-x$ are independent;
namely, $\pi(\-x)$ can be written as, 
 $$\pi(\-x)  = \prod_{i=1}^{d} \phi_i(x_i).$$  
We use the following algorithm to generate another sample $\-x^*$ from the MMC biasing distribution $f(\cdot)$. 
\begin{enumerate}
\item For $i = 1, \cdots, d$, sample $\xi_i \sim q_i(\cdot | x_i)$, where $q_i(\cdot)$ is a univariate PDF for $\xi_i$ centered at $x_i$ with the symmetry property $q_i(\xi_i | x_i) = q_k(x_i | \xi_i)$. 

\item Compute the acceptance probability $r_i = \min\{1,\, {\phi_i(\xi_i)}/{\phi_i(x_i)}\}$ for $i = 1, \cdots, d$, and then determine the $i$-th coordinate of the candidate sample by accepting or rejecting $\xi_i$ according to, 
\begin{equation}
\zeta_i = \left\{ 
  \begin{array} {rl}
    \xi_i, & \text{with probability }   r_i; \\
      x_i, & \text{with probability }  1 -r_i.
  \end{array}
\right.
\label{mcmc-zeta}
\end{equation}

\item Compute the final acceptance probability $r^* =\min\{ {{1,\,\Theta(\-x)}/\Theta(\bm{\zeta})}\}$,  and accept or reject the possible sample $\bm\zeta$ 
according to
\begin{equation}
\-x^* = \left\{ 
  \begin{array} {rl}
    \bm\zeta, & \text{with probability }  r^*; \\
      \-x, & \text{with probability }  1-r^*.
  \end{array}
\right.
\label{mcmc-new-x}
\end{equation}
\end{enumerate}
The ergodicity of the modified MCMC algorithm can be proved using the same arguments of \cite{AuB01} and so is omitted here. 
\bibliographystyle{plain}
\bibliography{mmc,reliability,rbo}

\end{document}